\documentclass[12pt]{amsart}

\usepackage[mathscr]{euscript}
\usepackage{mathtools}
\usepackage{amssymb}
\usepackage{stackrel}
\usepackage{pb-diagram,pb-xy} 
\usepackage[all,cmtip,arrow]{xy} 

\usepackage{quiver}
\usepackage{tikz-cd}

\usepackage[margin=1.5in]{geometry}  
\usepackage{hyperref} 
\hypersetup{colorlinks=true,linkcolor=blue,citecolor=magenta}

\numberwithin{equation}{section}

\theoremstyle{plain}
\newtheorem{theorem}[equation]{Theorem}
\newtheorem{lemma}[equation]{Lemma}
\newtheorem{corollary}[equation]{Corollary}
\newtheorem{proposition}[equation]{Proposition}

\theoremstyle{definition}
\newtheorem{definition}[equation]{Definition}

\newtheorem{remark}[equation]{Remark}

\newtheorem{example}[equation]{Example}

\newcommand{\homeq}{\simeq}

\newcommand{\N}{\mathbb{N}}

\newcommand{\weq}{\; \tilde{\longrightarrow} \;}

\mathchardef\mhyphen="2D

\newcommand{\colim}{\operatorname{colim} }

\newcommand{\Fun}{\operatorname{Fun} }
\newcommand{\End}{\operatorname{End} }

\newcommand{\Hom}{\operatorname{Hom} }
\newcommand{\Map}{\operatorname{Map} }

\newcommand{\Exc}{\operatorname{Exc} }

\newcommand{\un}[1]{\underline{#1}}

\newcommand{\cat}[1]{\mathscr{#1}}
\newcommand{\bcat}[1]{\mathbb{#1}}

\newcommand{\Weil}{\mathbb{W}\mathrm{eil}}

\newcommand{\spectra}{{\mathscr{S}p}}
\newcommand{\finbased}{\mathsf{S}_{\mathsf{fin,*}}}
\newcommand{\spaces}{\mathscr{S}}

\newcommand{\Catinf}{\mathbb{C}\mathrm{at}_\infty}
\newcommand{\Catdiff}{{\mathbb{C}\mathrm{at}_\infty^{\mathrm{diff}}}}

\newcommand{\Catpr}{\mathbb{C}\mathrm{at}_\infty^{\mathrm{pr}}}

\newcommand{\Topos}{\mathbb{T}\mathrm{opos}_\infty}
\newcommand{\InjTopos}{\mathbb{I}\mathrm{nj}\Topos}
\newcommand{\Logos}{\mathbb{L}\mathrm{ogos}_{\infty}}
\newcommand{\Catprca}{\Catinf^{\mathrm{pr,ca}}}
\newcommand{\Catacc}{\Catinf^{\mathrm{acc,\omega}}}

\newcommand{\Inf}{\mathbb{I}\mathrm{nf}}

\newcommand{\Sh}{\cat{S}h}
\newcommand{\PS}[1]{\cat{P}(\mathsf{#1})}

\newcommand{\tinfty}{\texorpdfstring{$\infty$}{infinity}}

\begin{document}
	
\title{Dual tangent structures for \tinfty-toposes}
	
\author{Michael Ching}
\email{mching@amherst.edu}
\address{Department of Mathematics and Statistics, Amherst College, PO Box 5000, Amherst, MA 01002, USA}
	
\begin{abstract}
We describe dual notions of tangent bundle for an $\infty$-topos, each underlying a tangent $\infty$-category in the sense of Bauer, Burke and the author. One of those notions is Lurie's tangent bundle functor for presentable $\infty$-categories, and the other is its adjoint. We calculate that adjoint for injective $\infty$-toposes, where it is given by applying Lurie's tangent bundle on $\infty$-categories of points.
\end{abstract}

\maketitle

In~\cite{bauer/burke/ching:2021}, Bauer, Burke, and the author introduce a notion of tangent $\infty$-category which generalizes the Rosick\'{y}~\cite{rosicky:1984} and Cockett-Cruttwell~\cite{cockett/cruttwell:2014} axiomatization of the tangent bundle functor on the category of smooth manifolds and smooth maps. We also constructed there a significant example: the Goodwillie tangent structure on the $\infty$-category $\Catdiff$ of (differentiable) $\infty$-categories, which is built on Lurie's tangent bundle functor. That tangent structure encodes the ideas of Goodwillie's calculus of functors~\cite{goodwillie:2003} and highlights the analogy between that theory and the ordinary differential calculus of smooth manifolds.

The goal of this note is to introduce two further examples of tangent $\infty$-categories: one on the $\infty$-category of $\infty$-toposes and geometric morphisms, which we denote $\Topos$, and one on the opposite $\infty$-category $\Topos^{op}$ which, following Anel and Joyal~\cite{anel/joyal:2019}, we denote $\Logos$.

These two tangent structures each encode a notion of tangent bundle for $\infty$-toposes, but from dual perspectives. As described in~\cite[Sec. 4]{anel/joyal:2019}, one can view $\infty$-toposes either from a `geometric' or `algebraic' point of view. In the former, an $\infty$-topos is thought of as a \emph{generalized topological space}; for example, each actual topological space gives rise to an $\infty$-topos $\Sh(X)$ of sheaves on $X$ (with values in the $\infty$-category $\spaces$ of $\infty$-groupoids), and each continuous map $X \to Y$ determines a geometric morphism $\Sh(X) \to \Sh(Y)$. From the algebraic perspective, an $\infty$-topos is more like a \emph{category of (higher) groupoids} with the $\infty$-topos $\spaces$ being a prime example. The `algebraic' morphisms between two $\infty$-toposes are those that preserve colimits and finite limits; i.e. the left adjoints of the geometric morphisms.

Our tangent structure on the `algebraic' $\infty$-category $\Logos$ is simply the restriction of the Goodwillie tangent structure. For an $\infty$-topos $\cat{X}$, the tangent bundle $T(\cat{X})$ is described by Lurie in~\cite[7.3.1]{lurie:2017}; the fibre of that bundle over each object $C \in \cat{X}$, i.e. the tangent space $T_C\cat{X}$, is the stabilization $\spectra(\cat{X}_{/C})$ of the slice $\infty$-topos over $C$.

The tangent structure on the `geometric' category $\Topos$ is \emph{dual} to that on $\Logos \homeq \Topos^{op}$ in a sense described by Cockett and Cruttwell~\cite[5.17]{cockett/cruttwell:2014}. Lurie's tangent bundle functor $T: \Topos^{op} \to \Topos^{op}$ has a left adjoint $U^{op}$ whose opposite $U: \Topos \to \Topos$ is the underlying functor for a tangent structure. We call that structure the \emph{geometric tangent structure} on $\Topos$.

The geometric tangent structure is \emph{representable} in the sense of~\cite[Sec. 5.2]{cockett/cruttwell:2014}. That is to say that the tangent bundle functor $U: \Topos \to \Topos$ is given by the exponential objects
\[ U(\cat{X}) = \cat{X}^{\cat{T}} \]
for some object $\cat{T}$ in the $\infty$-category $\Topos$, with the tangent structure on $U$ arising from so-called \emph{infinitesimal} structure on $\cat{T}$. This picture follows the same pattern as the tangent category associated to a model of Synthetic Differential Geometry (SDG)~\cite[5.10]{cockett/cruttwell:2014}, and we wonder which other features of SDG have a counterpart in the tangent $\infty$-category $\Topos$. That question is not explored here.

The infinitesimal object $\cat{T}$ is the $\infty$-topos of `parameterized spectra', also known as the Goodwillie tangent bundle $T(\spaces)$ on the $\infty$-category of spaces $\spaces$, with infinitesimal structure determined by the Goodwillie tangent structure. The exponential objects $U(\cat{X}) = \cat{X}^{T(\spaces)}$, however, do not seem to have a simple description, and we do not have a good understanding of the geometric tangent structure in its entirety.

We do offer a perspective on the geometric tangent structure by looking at its restriction to the subcategory of $\Topos$ consisting of those $\infty$-toposes that are \emph{injective} in a sense that generalizes Johnstone's notion of injective $1$-topos~\cite{johnstone:1981}. We prove that the geometric tangent structure on injective $\infty$-toposes is equivalent, via the `$\infty$-category of points' functor, to the Goodwillie tangent structure (restricted to those $\infty$-categories that are presentable and compactly-assembled~\cite[21.1.2]{lurie:2018}). We therefore view the geometric tangent structure as an \emph{extension} of the Goodwillie structure, to non-injective $\infty$-toposes, in addition to being its dual.

Here is an outline of the paper. In Section~\ref{sec:rep} we introduce the notion of `infinitesimal object' in an $\infty$-category $\bcat{X}$ with finite products (or, more generally, in any monoidal $\infty$-category), and we describe what it means for a tangent structure to be represented, or corepresented, by an infinitesimal object. These definitions extend to $\infty$-categories the corresponding notions for tangent categories due to Cockett and Cruttwell~\cite[Sec. 5]{cockett/cruttwell:2014}.

In Section~\ref{sec:top} we turn to $\infty$-toposes, and introduce the two perspectives given by the $\infty$-categories $\Topos$ and $\Logos = \Topos^{op}$. We then construct an infinitesimal object $\cat{T}$  in $\Topos$ whose underlying $\infty$-topos is $T(\spaces)$, and we show that $\cat{T}$ represents and corepresents tangent structures on $\Topos$ and $\Logos$ respectively. As mentioned, we prove that the latter of these is the restriction of the Goodwillie tangent structure of~\cite{bauer/burke/ching:2021} to the $\infty$-toposes.

Finally in Section~\ref{sec:tan} we give our description (Theorem~\ref{thm:Up}) of the geometric tangent structure on the subcategory of $\Topos$ consisting of the injective $\infty$-toposes, and make some explicit calculations that arise from that result.

\subsection*{Acknowledgements}

This paper is an extension of~\cite{bauer/burke/ching:2021}, and hence owes much to conversations with Kristine Bauer and Matthew Burke. The question of how the Goodwillie tangent structure is related to $\infty$-toposes was suggested by communications with Alexander Oldenziel and with Eric Finster. This work was supported by the National Science Foundation under grant DMS-1709032.


\section{Representable tangent \tinfty-categories} \label{sec:rep}

The goal of this section is to extend Cockett and Cruttwell's notion of \emph{representable} tangent category~\cite[Sec. 5.2]{cockett/cruttwell:2014} to the tangent $\infty$-categories of \cite{bauer/burke/ching:2021}. We start by developing the notion of \emph{infinitesimal object} in the sense of~\cite[5.6]{cockett/cruttwell:2014}. As with extending other tangent category notions to $\infty$-categories, we do this by re-expressing the Cockett-Cruttwell definition in terms of the category of Weil-algebras that Leung used in~\cite{leung:2017} to characterize tangent structures.

Throughout this paper we rely on notation and definitions from~\cite{bauer/burke/ching:2021}. In particular, let $\Weil$ be the monoidal category of \emph{Weil-algebras} of~\cite[1.1]{bauer/burke/ching:2021}. The objects of $\Weil$ are the augmented commutative semi-rings of the form
\[ \N[x_1,\dots,x_n]/(x_ix_j \; | \; i \sim j) \]
where the relations are quadratic monomials determined by equivalence relations on the sets $\{1,\dots,n\}$, for $n \geq 0$. Morphisms in $\Weil$ are semi-ring homomorphisms that commute with the augmentations, and monoidal structure is given by the tensor product which is also the coproduct. Certain pullback squares in the category $\Weil$ play a crucial role in the definition of tangent structure; see~\cite[1.10, 1.12]{bauer/burke/ching:2021}. We refer to those diagrams as the \emph{tangent pullbacks}.

We now give our notion of infinitesimal object in a monoidal $\infty$-category.

\begin{definition} \label{def:infinitesimal}
Let $\bcat{X}^{\boxtimes}$ be a monoidal $\infty$-category, and let $\bcat{X}^{op,\boxtimes}$ denote the corresponding monoidal structure on the opposite $\infty$-category $\bcat{X}^{op}$. An \emph{infinitesimal object} in $\bcat{X}^{\boxtimes}$ is a monoidal functor
\[ \cat{D}^\bullet: \Weil^{\otimes} \to \bcat{X}^{op,\boxtimes}.  \]
for which the underlying functor $\Weil \to \bcat{X}^{op}$ preserves the tangent pullbacks (i.e. maps tangent pullbacks in $\Weil$ to pushouts in $\bcat{X}$). These monoidal functors (and their monoidal natural transformations) form an $\infty$-category (see~\cite[2.1.3.7]{lurie:2017}), whose opposite\footnote{We take opposites here to ensure that morphisms in $\Inf(\bcat{X}^{\otimes})$ are built from morphisms in $\bcat{X}$ rather than $\bcat{X}^{op}$.}  we refer to as the \emph{$\infty$-category of infinitesimal objects in $\bcat{X}^{\boxtimes}$}, denoted $\Inf(\bcat{X}^{\boxtimes})$.
\end{definition}

\begin{example} \label{ex:inf-obj}
If $\bcat{X}$ is any $\infty$-category with finite products, then we refer to an infinitesimal object in the cartesian monoidal $\infty$-category $\bcat{X}^{\times}$ simply as an \emph{infinitesimal object in $\bcat{X}$}. We write $\Inf(\bcat{X})$ for the $\infty$-category of these infinitesimal objects.
\end{example}

\begin{example} \label{ex:inf-ord}
Let $\bcat{X}$ be an ordinary category with finite products. Our notion of infinitesimal object in $\bcat{X}$ agrees with that given by Cockett-Cruttwell~\cite[5.6]{cockett/cruttwell:2014} except that we do \emph{not} require that the objects $\cat{D}^A$ are exponentiable (the axiom there labelled [\textbf{Infsml.6}]). That condition is added in Proposition~\ref{prop:represented} below to explain when $\cat{D}^\bullet $ represents a tangent structure on $\bcat{X}$.
\end{example}

\begin{remark}
In the language of~\cite[5.7]{bauer/burke/ching:2021}, an infinitesimal object in a monoidal $\infty$-category $\bcat{X}^{\boxtimes}$ is precisely a tangent object in the $(\infty,2)$-category $\un{\bcat{X}}^{op,\boxtimes}$ that has a single object, mapping $\infty$-category $\bcat{X}^{op}$, and composition given by the monoidal structure $\boxtimes$.
\end{remark}

In the case that the monoidal structure on $\bcat{X}$ is given by the cartesian product, as in Example~\ref{ex:inf-obj}, there is a particularly simple way to identify infinitesimal objects. 

\begin{proposition} \label{prop:inf}
Let $\bcat{X}$ be an $\infty$-category with finite products. Taking underlying functors determines an equivalence between $\Inf(\bcat{X})$ and the full subcategory of $\Fun(\Weil^{op},\bcat{X})$ whose objects are the functors $\cat{D}^\bullet: \Weil^{op} \to \bcat{X}$ such that
\begin{enumerate}
	\item $\cat{D}^{\N}$ is a terminal object in $\bcat{X}$;
	\item for Weil-algebras $A,A'$ there is an equivalence in $\bcat{X}$
	\[ \cat{D}^{A \otimes A'} \weq \cat{D}^A \times \cat{D}^{A'}, \]
	induced by the canonical maps from $A$ and $A'$ to their coproduct $A \otimes A'$;
	\item $\cat{D}^\bullet$ maps the tangent pullbacks in $\Weil$ to pushouts in $\bcat{X}$.
\end{enumerate}
\end{proposition}
\begin{proof}
Since the monoidal structures on both $\Weil$ and $\bcat{X}^{op}$ are given by the coproduct, we can apply~\cite[2.4.3.8]{lurie:2017} to identify the $\infty$-category of \emph{lax} monoidal functors $\Weil^{\otimes} \to \bcat{X}^{op,\times}$ with the $\infty$-category of \emph{all} functors $\cat{D}^\bullet: \Weil \to \bcat{X}^{op}$. Conditions (1) and (2) on $\cat{D}^\bullet$ correspond to the case where that lax monoidal functor is monoidal, and (3) is the condition that this monoidal functor is an infinitesimal structure.
\end{proof}

We now consider two ways in which an infinitesimal object can determine a tangent structure, corresponding to those described in~\cite[Prop. 5.7 and Cor. 5.18]{cockett/cruttwell:2014}.

\begin{proposition} \label{prop:corepresented}
Let $\bcat{X}^{\boxtimes}$ be a monoidal $\infty$-category for which the monoidal structure $\boxtimes$ commutes with pushouts in its first variable, and let $\cat{D}^\bullet$ be an infinitesimal object in $\bcat{X}^{\boxtimes}$. Then there is a tangent structure on the $\infty$-category $\bcat{X}^{op}$ given by the $\Weil$-action map
\[ T: \Weil \times \bcat{X}^{op} \to \bcat{X}^{op}; \quad T^A(\cat{C}) = \cat{D}^A \boxtimes \cat{C}  \]
for a Weil-algebra $A$ and $\cat{C} \in \bcat{X}$.
\end{proposition}
\begin{proof}
The monoidal structure on $\bcat{X}$ determines a monoidal functor
\[ \boxtimes:  \bcat{X}^{op,\boxtimes} \to \End(\bcat{X}^{op})^{\circ}; \quad \cat{D} \mapsto \cat{D} \boxtimes - \]
which, by hypothesis, preserves pullbacks. Here $\End(\bcat{X}^{op})^{\circ}$ denotes the $\infty$-category of endofunctors on $\bcat{X}^{op}$ with monoidal structure given by composition.

Composing the infinitesimal object $\cat{D}^\bullet: \Weil^{\otimes} \to \bcat{X}^{op,\boxtimes}$ with the monoidal functor $\boxtimes$, we get a monoidal functor
\[ \cat{D}^\bullet \boxtimes - : \Weil^{\otimes} \to \End(\bcat{X}^{op})^{\circ} \]
which preserves the tangent pullbacks. Thus $T$ is a tangent structure on $\bcat{X}^{op}$; see~\cite[2.1]{bauer/burke/ching:2021}.
\end{proof}

\begin{definition} \label{def:corepresentable}
We say the tangent structure $T$ in Proposition~\ref{prop:corepresented} is \emph{corepresented} by the infinitesimal object $\cat{D}^\bullet$. A tangent structure is \emph{corepresentable} if it is equivalent to a tangent structure corepresented by some infinitesimal object.
\end{definition}

\begin{proposition} \label{prop:represented}
Let $\bcat{X}^{\boxtimes}$ be a monoidal $\infty$-category such that the monoidal product $\boxtimes$ preserves pushouts in its first variable. Let $\cat{D}^\bullet$ be an infinitesimal object in $\bcat{X}^{\boxtimes}$ such that for each Weil-algebra $A$, the functor $\cat{D}^A \boxtimes -: \bcat{X} \to \bcat{X}$ admits a right adjoint $\Map^{\boxtimes}_{\bcat{X}}(\cat{D}^A,-)$. Then there is a tangent structure on the $\infty$-category $\bcat{X}$ given by
\[ U: \Weil \times \bcat{X} \to \bcat{X}; \quad U^A(\cat{C}) = \Map^{\boxtimes}_{\bcat{X}}(\cat{D}^A,\cat{C}). \]
\end{proposition}
\begin{proof}
Note that we do not assume that the monoidal structure $\boxtimes$ as a whole is \emph{closed}, only that certain specific functors admit a right adjoint. Our first task is to show that those right adjoints can be chosen functorially.

Let $\End^L(\bcat{X})^{\circ} \subseteq \End(\bcat{X})^{\circ}$ be the full (monoidal) subcategory whose objects are the left adjoint functors $\bcat{X} \to \bcat{X}$, and similarly for $\End^R(\bcat{X})$. Then there is an equivalence of monoidal $\infty$-categories
\[ \mathrm{adj}: \End^L(\bcat{X})^{\circ} \weq \End^R(\bcat{X})^{\circ} \]
that sends a functor to some choice of its right adjoint. Such an equivalence can be constructed in a similar manner to that in~\cite[5.5.3.4]{lurie:2009}.

By hypothesis, the composite
\[ \dgTEXTARROWLENGTH=2.5em \Weil^{\otimes} \arrow{e,t}{\cat{D}^\bullet} \bcat{X}^{op,\boxtimes} \arrow{e,t}{\boxtimes}  \End(\bcat{X}^{op})^{\circ} \arrow{e,tb}{\sim}{op} \End(\bcat{X})^{\circ} \]
takes values in $\End^L(\bcat{X})^{\circ}$. Therefore we can form the composite monoidal functor
\[ \dgTEXTARROWLENGTH=2.5em \Weil^{\otimes} \to \End^L(\bcat{X})^{\circ} \arrow{e,tb}{\sim}{\mathrm{adj}} \End^R(\bcat{X})^{\circ} \to \End(\bcat{X})^{\circ}.  \]
It remains to show that the underlying functor $A \mapsto \Map_{\bcat{X}}^{\boxtimes}(\cat{D}^A,-)$ preserves each of the tangent pullbacks in $\Weil$. Suppose
\[\begin{tikzcd}
	A & {A_1} \\
	{A_2} & {A_0}
	\arrow[from=1-1, to=1-2]
	\arrow[from=1-1, to=2-1]
	\arrow[from=2-1, to=2-2]
	\arrow[from=1-2, to=2-2]
\end{tikzcd}\]
is one of those tangent pullbacks. Then it is sufficient to show that for all $\cat{C},\cat{E} \in \bcat{X}$, the following diagram is a pullback of mapping spaces:
\[\begin{tikzcd}
	\Hom_{\bcat{X}}(\cat{E},\Map^{\boxtimes}_{\bcat{X}}(\cat{D}^A,\cat{C})) &
	{\Hom_{\bcat{X}}(\cat{E},\Map^{\boxtimes}_{\bcat{X}}(\cat{D}^{A_1},\cat{C}))} \\
	{\Hom_{\bcat{X}}(\cat{E},\Map^{\boxtimes}_{\bcat{X}}(\cat{D}^{A_2},\cat{C}))} & {\Hom_{\bcat{X}}(\cat{E},\Map^{\boxtimes}_{\bcat{X}}(\cat{D}^{A_0},\cat{C}))}
	\arrow[from=1-1, to=1-2]
	\arrow[from=1-1, to=2-1]
	\arrow[from=2-1, to=2-2]
	\arrow[from=1-2, to=2-2]
\end{tikzcd}\]
We can write this diagram equivalently as
\[\begin{tikzcd}
	\Hom_{\bcat{X}}(\cat{D}^{A} \boxtimes \cat{E}, \cat{C}) & \Hom_{\bcat{X}}(\cat{D}^{A_1} \boxtimes \cat{E}, \cat{C}) \\
	\Hom_{\bcat{X}}(\cat{D}^{A_2} \boxtimes \cat{E}, \cat{C}) & \Hom_{\bcat{X}}(\cat{D}^{A_0} \boxtimes \cat{E}, \cat{C})
	\arrow[from=1-1, to=1-2]
	\arrow[from=1-1, to=2-1]
	\arrow[from=2-1, to=2-2]
	\arrow[from=1-2, to=2-2]
\end{tikzcd}\]
which is a pullback since $- \boxtimes \cat{E}$ preserves pushouts by hypothesis.
\end{proof}

\begin{definition} \label{def:representable}
We say that the tangent structure $U$ in Proposition~\ref{prop:represented} is \emph{represented} by the infinitesimal object $\cat{D}^\bullet$. A tangent structure is \emph{representable} if it is equivalent to one represented by some infinitesimal object.
\end{definition}

\begin{definition} \label{def:dual}
Tangent structures on $\bcat{X}$ and $\bcat{X}^{op}$ are \emph{dual} if they are, respectively, represented and corepresented by the same infinitesimal object. It follows from comparing the hypotheses in Propositions~\ref{prop:corepresented} and~\ref{prop:represented} that any representable tangent structure on $\bcat{X}$ has a dual tangent structure on $\bcat{X}^{op}$.
\end{definition}

\section{\tinfty-toposes} \label{sec:top}

The main purpose of this paper is to construct dual tangent structures on a certain $\infty$-category $\Topos$ and its opposite, whose objects are $\infty$-toposes. In this section we introduce that $\infty$-category and construct the infinitesimal object $\cat{T}^\bullet$ that (co)represents those tangent structures. Our main reference for the $\infty$-category $\Topos$ is~\cite[Sec. 6.3]{lurie:2009} where it is denoted $\mathscr{RT}\mathrm{op}_\infty$.

To define the $\infty$-category $\Topos$ we have to pay some attention to size issues. We assume three nested Grothendieck universes and refer to simplicial sets of these sizes as \emph{small}, \emph{large} and \emph{very large}, respectively. Let $\Catinf$ denote the (very large) $\infty$-category of large $\infty$-categories. One of the objects in $\Catinf$ is the (large) $\infty$-category $\spaces$ of small spaces, i.e. Kan complexes.

\begin{definition} \label{def:topos}
An \emph{$\infty$-topos} is a (large) $\infty$-category $\cat{X}$ that is an accessible left exact localization of the $\infty$-category $\PS{C} = \Fun(\mathsf{C}^{op},\spaces)$ of presheaves on some small $\infty$-category $\mathsf{C}$. In other words, there is some small $\infty$-category $\mathsf{C}$ and an adjunction
\begin{equation} \label{eq:adj} \begin{tikzcd}
	{\PS{C}} & {\cat{X}}
	\arrow[""{name=0, anchor=center, inner sep=0}, "f", shift left=2, from=1-1, to=1-2]
	\arrow[""{name=1, anchor=center, inner sep=0}, "g", shift left=2, from=1-2, to=1-1]
	\arrow["\dashv"{anchor=center, rotate=-90}, draw=none, from=0, to=1]
\end{tikzcd} \end{equation}
such that $g$ is fully faithful and accessible (preserves  $\kappa$-filtered colimits for some small regular cardinal $\kappa$), and $f$ preserves finite limits.
\end{definition}

\begin{example}
The presheaf $\infty$-category $\cat{P}(\mathsf{C})$, for a small $\infty$-category $\mathsf{C}$, is an $\infty$-topos. In particular $\spaces$ is an $\infty$-topos which, by~\cite[6.3.4.1]{lurie:2009}, is a terminal object in the $\infty$-category $\Topos$ which we now introduce.
\end{example}

\begin{definition}
Let $\Topos$ denote the subcategory of $\Catinf$ whose objects are the $\infty$-toposes and whose morphism are the \emph{geometric morphisms}, i.e. those functors $F: \cat{X} \to \cat{Y}$ which admit a left adjoint $\cat{Y} \to \cat{X}$ that preserves finite limits.
\end{definition}

The opposite $\infty$-category $\Topos^{op}$ is also equivalent to a subcategory of $\Catinf$ which, following Anel and Joyal~\cite{anel/joyal:2019}, we denote by $\Logos$.

\begin{definition}
Let $\Logos$ be the subcategory of $\Catinf$ whose objects are the $\infty$-toposes and morphisms are the functors that preserve small colimits and finite limits. There is an equivalence $\Logos \homeq \Topos^{op}$ that is the identity on objects and maps a functor $\cat{Y} \to \cat{X}$ in $\Logos$ to the right adjoint guaranteed by the Adjoint Functor Theorem~\cite[5.5.2.9]{lurie:2009}, a geometric morphism $\cat{X} \to \cat{Y}$; see~\cite[6.3.1.8]{lurie:2009}.
\end{definition}

The remainder of this section is devoted to the construction of an infinitesimal object in $\Topos \homeq \Logos^{op}$, i.e. a monoidal functor
\[ \cat{T}^\bullet: \Weil^{\otimes} \to \Logos^{\boxtimes} \]
where $\boxtimes$ denotes the coproduct in the $\infty$-category $\Logos$, or equivalently the product in $\Topos$.

The infinitesimal object $\cat{T}^\bullet$ is derived from the Goodwillie tangent structure constructed in~\cite{bauer/burke/ching:2021} which we now recall. See~\cite[2.1]{bauer/burke/ching:2021} for the notion of tangent structure on an $\infty$-category, and see~\cite[Sec. 7]{bauer/burke/ching:2021} for the basic facets of the Goodwillie tangent structure.

\begin{definition} \label{def:goodwillie}
Let $\Catdiff$ be the $\infty$-category of (large) differentiable\footnote{An $\infty$-category is \emph{differentiable} if it has finite limits and sequential colimits which commute. Any $\infty$-topos is differentiable by~\cite[7.3.4.7]{lurie:2009}.} $\infty$-categories and sequential-colimit-preserving functors. Let $\finbased$ denote the (small) $\infty$-category of pointed finite spaces.

The \emph{Goodwillie tangent structure} on $\Catdiff$ is a map
\[ T: \Weil \times \Catdiff \to \Catdiff \]
given on a Weil-algebra $A$ with $n$ generators, and differentiable $\infty$-category $\cat{C}$, by the subcategory
\begin{equation} \label{eq:PA} T^A(\cat{C}) = \Exc^A(\finbased^n,\cat{C}) \subseteq \Fun(\finbased^n,\cat{C}) \end{equation}
of functors $\finbased^n \to \cat{C}$ that are \emph{$A$-excisive}\footnote{A functor is \emph{excisive} if it maps pushout squares to pullbacks. The notion of $A$-excisive is a multivariable generalization of excisive that reflects the structure of the Weil-algebra $A$.} in the sense described in~\cite[7.1]{bauer/burke/ching:2021}. By~\cite[7.5]{bauer/burke/ching:2021}, the inclusion (\ref{eq:PA}) admits a left adjoint $P_A$ which preserves finite limits.

The action of $T$ on morphisms (in $\Weil$ and $\Catdiff$) is described in detail in~\cite[7.7 and 7.14]{bauer/burke/ching:2021}. For a Weil-algebra morphism $\phi: A \to A'$ and (sequential-colimit-preserving) functor $G: \cat{C} \to \cat{D}$, we have
\[ T^{\phi}(F): T^A(\cat{C}) \to T^{A'}(\cat{D}); \quad L \mapsto P_{A'}(GL\tilde{\phi}) \]
where $\tilde{\phi}: \finbased^{n'} \to \finbased^n$ is a functor built to the same pattern as the algebra homomorphism $\phi$; see~\cite[7.12]{bauer/burke/ching:2021}.
\end{definition}

\begin{proposition} \label{prop:Goodwillie-Logos}
The Goodwillie tangent structure on $\Catdiff$
 restricts to a tangent structure on the subcategory $\Logos \subseteq \Catdiff$.
\end{proposition}
\begin{proof}
Suppose first that $\cat{X}$ is an $\infty$-topos and $A$ is a Weil-algebra with $n$ generators. By~\cite[7.5]{bauer/burke/ching:2021}, the $\infty$-category $T^A(\cat{X})$ is an accessible left exact localization of the $\infty$-topos $\Fun(\finbased^n,\cat{X})$, hence $T^A(\cat{X})$ is an $\infty$-topos.

Now let $G: \cat{X} \to \cat{Y}$ be a morphism in $\Logos$ and $\phi:A \to A'$ a Weil-algebra morphism. We have to show that the functor
\[ T^\phi(G) : T^A(\cat{X}) \to T^{A'}(\cat{Y}); \quad L \mapsto P_{A'}(GL\tilde{\phi}) \]
is also in $\Logos$. Finite limits in $T^A(\cat{X})$ and $T^{A'}(\cat{Y})$ are calculated objectwise, and both $G$ and $P_{A'}$ preserve those finite limits, so $T^{\phi}(G)$ preserves finite limits.
	
Let $(L_{\alpha})$ be a diagram in $T^A(\cat{X})$ with colimit $L$. Then we have an equivalence $L \homeq P_A(\colim L_{\alpha})$ where $\colim$ denotes the (objectwise) colimit calculated in $\Fun(\finbased^n,\cat{C})$. Then there is a sequence of equivalences:
\[ \begin{split} P_{A'}(GL\tilde{\phi}) &\homeq P_{A'}(GP_A(\colim L_{\alpha})\tilde{\phi}) \\ 
		&\homeq P_{A'}(G(\colim L_{\alpha})\tilde{\phi}) \\
		&\homeq P_{A'}(\colim GL_{\alpha}\tilde{\phi}) \end{split} \]
where we have used equivalences of the form~\cite[7.11 and 7.25]{bauer/burke/ching:2021} to identify the first and second lines, and the fact that $G$ preserves colimits to identify the second and third.

Therefore $T^\phi(G)$ preserves colimits, which completes the proof that $T^{\phi}(G)$ is a morphism in $\Logos$, and hence that the $\Weil$-action on $\Catdiff$ restricts to a functor
\[ T: \Weil \times \Logos \to \Logos. \]	
It remains to show that $T$ preserves the tangent pullbacks in $\Weil$. Since limits in $\Logos$ are calculated in $\Catinf$ by~\cite[6.3.2.3]{lurie:2009}, that claim follows from~\cite[7.36 and 7.38]{bauer/burke/ching:2021}.
\end{proof}

We now construct the desired infinitesimal object.

\begin{definition} \label{def:T}
Define a functor
\[ \cat{T}^\bullet: \Weil \to \Logos \]
by
\[ \cat{T}^A := T^A(\spaces), \]
i.e. by evaluating the Goodwillie tangent structure at the $\infty$-topos $\spaces$ of spaces.
\end{definition}

\begin{proposition} \label{prop:T}
The functor $\cat{T}^\bullet$ is the underlying functor of an infinitesimal object in $\Topos \homeq \Logos^{op}$. 
\end{proposition}
\begin{proof}
We apply Proposition~\ref{prop:inf}:
\begin{enumerate} \itemsep=10pt
\item By~\cite[6.3.4.1]{lurie:2009}, $\cat{T}^{\N} = \spaces$ is the terminal object in $\Topos$.

\item The canonical map
\[ T^{A \otimes A'}(\spaces) \to T^A(\spaces) \boxtimes T^{A'}(\spaces) \]
is an equivalence of $\infty$-toposes; this claim follows from Lemma~\ref{lem:TAX} below by taking $\cat{X} = T^{A'}(\spaces)$ and noting that $T^{A \otimes A'}(\spaces) = T^A(T^{A'}(\spaces))$.

\item For each tangent pullback in $\Weil$, the corresponding diagram
\[\begin{tikzcd}
	{T^A(\spaces)} & {T^{A_1}(\spaces)} \\
	{T^{A_2}(\spaces)} & {T^{A_0}(\spaces)}
	\arrow[from=1-1, to=1-2]
	\arrow[from=1-1, to=2-1]
	\arrow[from=2-1, to=2-2]
	\arrow[from=1-2, to=2-2]
\end{tikzcd}\]
is a pullback in $\Logos$ (and hence a pushout in $\Topos$). This claim is part of the condition that $T$ is a tangent structure on $\Logos$, as verified in the last part of the proof of Proposition~\ref{prop:Goodwillie-Logos}.
\end{enumerate}
\end{proof}

To show that the infinitesimal object $\cat{T}^\bullet$ represents, and corepresents, tangent structures on the $\infty$-categories $\Topos$ and $\Logos$ respectively, we verify the conditions of Propositions~\ref{prop:corepresented} and~\ref{prop:represented}.

\begin{proposition} \label{prop:Topos-corep}
The product $\boxtimes$ on $\Topos$ preserves pushouts in each variable individually.
\end{proposition}
\begin{proof}
We use the fact, e.g. see~\cite[2.15]{anel/lejay:2018}, that the coproduct $\boxtimes$ in $\Logos$ is given by the tensor product of cocomplete $\infty$-categories; see~\cite[4.8.1]{lurie:2017}. Then~\cite[4.24]{anel/lejay:2018} tells us that for $\infty$-toposes $\cat{Y},\cat{X}$, we have
\[ \cat{Y} \boxtimes \cat{X} \homeq \Fun^{\lim}(\cat{X}^{op},\cat{Y}) \]
where the right-hand side is the $\infty$-category functors $\cat{Y}^{op} \to \cat{X}$ that preserve small limits.

We therefore have to show that $\Fun^{\lim}(\cat{X}^{op},-)$ preserves pushouts in $\Topos$ which, by~\cite[6.3.2.3]{lurie:2009}, are pullbacks in $\Catinf$. That claim is a consequence of~\cite[6.4.12]{riehl/verity:2020} which implies that pullbacks in the $\infty$-cosmos of $\infty$-categories with small limits are given by pullbacks in the $\infty$-cosmos of all $\infty$-categories.
\end{proof}

\begin{proposition} \label{prop:Topos-rep}
For each Weil-algebra $A$, the functor
\[ \cat{T}^A \boxtimes - : \Topos \to \Topos \]
admits a right adjoint.
\end{proposition}
\begin{proof}
Anel and Lejay~\cite[4.37]{anel/lejay:2018} show that any compactly-generated $\infty$-topos is exponentiable, so it is sufficient to show that $\cat{T}^A = \Exc^A(\finbased^n,\spaces)$ is compactly-generated. It follows from~\cite[5.3.5.12]{lurie:2009} that the presheaf $\infty$-category $\Fun(\finbased^n,\spaces)$ is compactly-generated, so by~\cite[5.5.7.3]{lurie:2009} it is sufficient to note that $\Exc^A(\finbased,\spaces)$ is closed under filtered colimits in $\Fun(\finbased^n,\spaces)$, which follows from the fact that filtered colimits in $\spaces$ commute with pullbacks.
\end{proof}

\begin{theorem} \label{thm:TS}
The infinitesimal object $\cat{T}^\bullet$ represents a tangent structure $U$ on the $\infty$-category $\Topos$ and corepresents a tangent structure $T$ on the $\infty$-category $\Logos$.
\end{theorem}
\begin{proof}
We apply Propositions~\ref{prop:represented} and~\ref{prop:corepresented}, respectively, using the results of Propositions~\ref{prop:Topos-corep} and~\ref{prop:Topos-rep}.
\end{proof}

We refer to the tangent structure $U$ on $\Topos$ as the \emph{geometric tangent structure}, and we begin the study of that structure in the next section. The corepresented tangent structure $T$ turns out to be much more familiar.

\begin{proposition}
The corepresentable tangent structure of Theorem~\ref{thm:TS} is equivalent to the restriction of the Goodwillie tangent structure of~\cite{bauer/burke/ching:2021} to the subcategory $\Logos \subseteq \Catdiff$, as described in Proposition~\ref{prop:Goodwillie-Logos}.
\end{proposition}
\begin{proof}
We have to show that the following diagram of monoidal functors commutes up to monoidal equivalence.
\[\begin{tikzcd}
	{\Weil^{\otimes}} & {\Logos^{\boxtimes}} \\
	& {\End(\Logos)^{\circ}}
	\arrow["{\cat{T}^\bullet}", from=1-1, to=1-2]
	\arrow["\boxtimes", from=1-2, to=2-2]
	\arrow["T"', from=1-1, to=2-2]
\end{tikzcd}\]
This claim is a consequence of the following lemma.
\end{proof}

\begin{lemma} \label{lem:TAX}
Let $\cat{X}$ be an $\infty$-topos, and $A$ a Weil-algebra. Then there is a canonical equivalence in $\Logos$
\[ T^A(\spaces) \boxtimes \cat{X} \weq T^A(\cat{X}) \]
built from the maps $T^A(!) : T^A(\spaces) \to T^A(\cat{X})$ and $T^{\eta}(\cat{X}): \cat{X} \to T^A(\cat{X})$, where $! : \spaces \to \cat{X}$ and $\eta: \N \to A$ are the maps from the initial objects in $\Logos$ and $\Weil$ respectively.
\end{lemma}
\begin{proof}
Using the approach from the proof of Proposition~\ref{prop:Topos-corep}, it is sufficient to show that the canonical map
\[ T^A(\cat{X}) \to \Fun^{\lim}(\cat{X}^{op},T^A(\spaces)); \quad L \mapsto \Hom_{\cat{X}}(-,L) \]
is an equivalence. Noting that limits in $T^A(\spaces) = \Exc^A(\finbased^n,\spaces)$ are calculated objectwise, we can rewrite the target of that map as the $\infty$-category of $A$-excisive functors $\finbased^n \to \Fun^{\lim}(\cat{X}^{op},\spaces)$, since pullbacks in $\Fun^{\lim}(\cat{X}^{op},\spaces)$ are also calculated objectwise. Our claim then follows by noting that the map
\[ \cat{X} \to \Fun^{\lim}(\cat{X}^{op},\spaces); \quad \Hom_{\cat{X}}(-,L) \]
is an equivalence by~\cite[4.24]{anel/lejay:2018} again.
\end{proof}

\section{The geometric tangent structure} \label{sec:tan}

The aim of this section is to begin a study of the geometric tangent structure on $\Topos$ given by Thorem~\ref{thm:TS}. By definition, the tangent bundle construction for this tangent structure,
\[ U(\cat{X}) = \cat{X}^{T(\spaces)}, \]
is the exponential object for the $\infty$-topos $T(\spaces)$. These exponential objects do not appear to be easy to calculate, though a general construction can be gleaned from the proofs of~\cite[4.33]{anel/lejay:2018} or~\cite[21.1.6.12]{lurie:2018}.

Those approaches proceed by first calculating the exponential object for a collection of $\infty$-toposes that are injective in the sense of Definition~\ref{def:inj} below. Since any $\infty$-topos $\cat{X}$ can be written as the pullback of a diagram of injective $\infty$-toposes~\cite[21.1.6.16]{lurie:2018}, and since $U$ preserves pullbacks, one might be able to recover an explicit description of $U(\cat{X})$ from those calculations, though we do not attempt that here.

It turns out that the geometric tangent structure for injective $\infty$-toposes has a compelling description. We prove in Theorem~\ref{thm:Up} that the `$\infty$-category of points' construction determines an equivalence between that geometric tangent structure and the Goodwillie tangent structure restricted to the presentable compactly-assembled $\infty$-categories of~\cite[21.1.2]{lurie:2018}. Thus, on injective $\infty$-toposes at least, one can view the geometric tangent structure as simply a different incarnation of the Goodwillie structure.

\begin{definition} \label{def:inj}
An $\infty$-topos $\cat{X}$ is \emph{injective} if $\cat{X}$ is a retract, in $\Topos$, of a presheaf $\infty$-category $\cat{P}(\mathsf{D})$ where $\mathsf{D}$ is a small $\infty$-category that has finite limits. Let $\InjTopos$ be the full subcategory of $\Topos$ consisting of the injective $\infty$-toposes.
\end{definition}

\begin{remark}
An $\infty$-topos is injective if and only if it satisfies the equivalent conditions of~\cite[21.1.5.4]{lurie:2018}; our definition is an intermediate step in proving (4) implies (1) in that result. Our definition is also equivalent to that of Anel and Lejay in~\cite[4.6]{anel/lejay:2018}; combine (4) of~\cite[21.1.5.4]{lurie:2018} with~\cite[2.6]{anel/lejay:2018}.
\end{remark}

The attraction of injective $\infty$-toposes is that they can be recovered from their $\infty$-categories of `points'.

\begin{definition}
The \emph{$\infty$-category of points} of an $\infty$-topos $\cat{X}$ is the $\infty$-category
\[ p(\cat{X}) := \Fun^*(\cat{X},\spaces) \]
of functors $\cat{X} \to \spaces$ that preserve small colimits and finite limits, i.e. the geometric morphisms $\spaces \to \cat{X}$. Since $\spaces$ is the terminal object in $\Topos$, the objects of $p(\cat{X})$ are indeed the `generalized points' of the $\infty$-topos $\cat{X}$. By~\cite[21.1.1.6]{lurie:2018}, the construction of $p(\cat{X})$ extends to a functor
\[ p: \Topos \to \Catacc \]
whose target is the subcategory of $\Catinf$ consisting of the $\infty$-categories that are accessible and admit filtered colimits, with morphisms the filtered-colimit-preserving functors.
\end{definition}

\begin{proposition} \label{prop:inj}
The functor $p$ restricts to an equivalence of $\infty$-categories
\[ p: \InjTopos \weq \Catprca \subseteq \Catacc \]
whose target consists of those $\infty$-categories $\cat{C}$ that are both presentable and compactly-assembled, in the sense of~\cite[21.1.2.1]{lurie:2018}. The inverse to $p$ maps such an $\infty$-category $\cat{C}$ to the $\infty$-topos $\Fun^{\omega}(\cat{C},\spaces)$ of filtered-colimit-preserving functors $\cat{C} \to \spaces$.
\end{proposition}
\begin{proof}
The inverse map $\Fun^{\omega}(-,\spaces)$ is fully faithful by~\cite[21.1.5.3]{lurie:2018}, essentially surjective by~\cite[21.1.5.4(1)]{lurie:2018}, and has inverse $p$ by~\cite[21.1.5.1]{lurie:2018}. 
\end{proof}

\begin{remark} \label{rem:catpresca}
An alternative approach to the proof of Proposition~\ref{prop:inj} is in~\cite[4.9]{anel/lejay:2018} which identifies $\Catprca$ with the full subcategory of $\Catacc$ consisting of retracts of the preasheaf $\infty$-categories. 
\end{remark}

We now prove the main result of this section, giving a calculation of the geometric tangent structure for injective $\infty$-toposes.

\begin{theorem} \label{thm:Up}
The equivalence $p$ of Proposition~\ref{prop:inj} underlies an equivalence of tangent structures
\[ p: (\InjTopos,U) \weq (\Catprca,T) \]
between the geometric tangent structure on $\InjTopos \subseteq \Topos$ and the Goodwillie tangent structure on $\Catprca$.\footnote{It is unclear to this author whether an arbitrary presentable compactly-assembled $\infty$-category $\cat{C}$ is differentiable, so $\Catprca$ is perhaps not a subcategory of $\Catdiff$. However, the construction of the Goodwillie tangent structure in~\cite{bauer/burke/ching:2021} can be carried out with $\Catdiff$ replaced by the $\infty$-category $\Catpr$ of presentable $\infty$-categories and filtered-colimit-preserving functors, of which $\Catprca$ is a full subcategory. Alternatively the reader may restrict attention to the compactly-generated $\infty$-categories which correspond to those $\infty$-toposes that are presheaves on a small $\infty$-category that has finite limits.}
\end{theorem}
\begin{proof}
We start by showing that the Goodwillie tangent structure restricts to $\Catprca$. Suppose $\cat{C}$ is presentable and compactly-assembled. Then, by Remark~\ref{rem:catpresca}, $\cat{C}$ is a retract, in $\Catacc$, of a presheaf $\infty$-category. Hence $\Fun(\finbased^n,\cat{C})$ is a retract of a presheaf $\infty$-category, so is also presentable and compactly-assembled. Finally, the map $P_A$ of Definition~\ref{def:goodwillie} displays $T^A(\cat{C})$ as a retract, in $\Catacc$, of $\Fun(\finbased^n,\cat{C})$, so $T^A(\cat{C})$ is also presentable and compactly-assembled.

Now let $q: \Catprca \to \InjTopos$ be the inverse to $p$ given by $q(\cat{C}) = \Fun^{\omega}(\cat{C},\spaces)$. We then define natural equivalences
\[ \alpha: qT^A \weq U^Aq \]
with components
\[ \alpha_{\cat{C}}: \Fun^{\omega}(T^A(\cat{C}),\spaces) \to U^A(\Fun^{\omega}(\cat{C},\spaces)) = \Fun^{\omega}(\cat{C},\spaces)^{T^A(\spaces)} \]
as follows. 

First note that the proof of Lemma~\ref{lem:TAX} relies purely on the identification of $\boxtimes$ with the tensor product for presentable $\infty$-categories, and so extends to give a canonical equivalence
\[ T^A(\spaces) \boxtimes \cat{C} \weq T^A(\cat{C}) \]
for all presentable $\infty$-categories. By~\cite[21.1.4.3]{lurie:2018}, and the argument of~\cite[21.1.6.9]{lurie:2018}, we also have equivalences of $\infty$-toposes of the form
\[ \Fun^{\omega}(\cat{X} \boxtimes \cat{C}, \spaces) \weq \Fun^{\omega}(\cat{C},\spaces)^{\cat{X}} . \]
Combining these two maps, with $\cat{X} = T^A(\spaces)$, yields the desired equivalence $\alpha_{\cat{C}}$. Note that the existence of these equivalences also verifies that $U$ restricts to a tangent structure on the subcategory $\InjTopos \subseteq \Topos$.

The construction of $\alpha_{\cat{C}}$ is natural (in $A$ and $\cat{C}$) and monoidal (with respect to the tensor product of Weil-algebras), so the maps $\alpha$ yield an equivalence of tangent structures with underlying functor $q$, whose inverse is the required tangent equivalence $p$.
\end{proof}

\begin{corollary} \label{cor:UX}
For an injective $\infty$-topos $\cat{X}$:
\[ U\cat{X} \homeq \Fun^{\omega}(T(p\cat{X}),\spaces). \]
\end{corollary}

\begin{corollary} \label{cor:UxX}
Let $\cat{X}$ be an injective $\infty$-topos, and let $x: \spaces \to \cat{X}$ be a generalized point in $\cat{X}$. Then the geometric tangent space $U_x\cat{X}$ (in $\Topos$) exists and has $\infty$-category of points
\[ p(U_x\cat{X}) \homeq T_x(p\cat{X}). \]
Note, however, that we have no reason to believe that $U_x\cat{X}$ is injective.
\end{corollary}
\begin{proof}
By definition the tangent space is the pullback in $\Topos$ of the form
\[\begin{tikzcd}
	{U_x\cat{X}} & {U\cat{X}} \\
	\spaces & {\cat{X}}
	\arrow[from=1-1, to=1-2]
	\arrow[from=1-1, to=2-1]
	\arrow["{\epsilon_{\cat{X}}}", from=1-2, to=2-2]
	\arrow["x", from=2-1, to=2-2]
\end{tikzcd}\]
This pullback exists, and is preserved by each $U^A$, since $\Topos$ has all limits, and $U^A$ is a right adjoint. The functor $p: \Topos \to \Catacc$ is a right adjoint by~\cite[21.1.1.6]{lurie:2018}, so applying $p$ we get a pullback diagram in $\Catacc$, and hence in $\Catinf$, of the form
\[\begin{tikzcd}
	{p(U_x\cat{X})} & {T(p\cat{X})} \\
	{*} & {p\cat{X}}
	\arrow[from=1-1, to=1-2]
	\arrow[from=1-1, to=2-1]
	\arrow["{\epsilon_{p\cat{X}}}", from=1-2, to=2-2]
	\arrow["x", from=2-1, to=2-2]
\end{tikzcd}\]
which identifies $p(U_x\cat{X})$ with the Goodwillie tangent space $T_x(p\cat{X})$, as claimed.
\end{proof}

Lurie's proof that a compactly-generated $\infty$-topos is exponentiable relies on a lemma~\cite[21.1.6.16]{lurie:2018} that says every $\infty$-topos is a pullback of injective $\infty$-toposes; see also~\cite[2.8]{anel/lejay:2018}. Since the geometric tangent bundle functor $U$ is a right adjoint, we could in principle use such pullbacks, together with the calculation in Corollary~\ref{cor:UX}, to give an explicit description of $U\cat{X}$ for any $\infty$-topos $\cat{X}$. 

Theorem~\ref{thm:Up} also gives us a different perspective on the relationship between the geometric and Goodwillie tangent structures. We saw in Theorem~\ref{thm:TS} that those tangent structures are dual, but now we can also view the geometric tangent structure on $\Topos$ as an \emph{extension} of the Goodwillie tangent structure. If we think of an $\infty$-topos as an $\infty$-category (of points) together with additional information, then the geometric tangent bundle on an injective $\infty$-topos simply \emph{is} the Goodwillie tangent bundle. The full geometric tangent structure on $\Topos$ extends the Goodwillie structure to $\infty$-toposes for which that additional information is nontrivial.

We conclude by giving some simple calculations of the geometric tangent structure based on Corollary~\ref{cor:UX}.

\begin{example}
For the terminal $\infty$-topos $\spaces$, we have
\[ U(\spaces) \homeq \spaces. \]
Since $U$ is a right adjoint and $\spaces$ is a terminal object, we did not need Theorem~\ref{thm:Up} to prove this fact. However, we now see that it corresponds to the calculation $T(*) \homeq *$ for the Goodwillie tangent bundle on the trivial $\infty$-category.
\end{example}

\begin{example} \label{ex:affine}
Let $\mathsf{C}$ be a small $\infty$-category, and let
\[ \cat{A}^{\mathsf{C}} := \cat{P}(\bar{\mathsf{C}}) = \Fun(\bar{\mathsf{C}}^{op},\spaces) \]
be the \emph{affine} $\infty$-topos of~\cite[2.7]{anel/lejay:2018}, where $\bar{\mathsf{C}}$ is obtained by freely adding finite limits to $\mathsf{C}$. The $\infty$-category of points of $\cat{A}^{\mathsf{C}}$ is
\[ p(\cat{A}^{\mathsf{C}}) \homeq \cat{P}(\mathsf{C}^{op}) = \Fun(\mathsf{C},\spaces) \]
and so the geometric tangent bundle is given by
\[ U(\cat{A}^{\mathsf{C}}) \homeq \Fun^{\omega}(\Fun(\mathsf{C},T(\spaces)),\spaces). \]
By Corollary~\ref{cor:UxX}, the geometric tangent space $U_x(\cat{A}^{\mathsf{C}})$, for a functor $x: \mathsf{C} \to \spaces$, has $\infty$-category of points
\[ p(U_x\cat{A}^{\mathsf{C}}) \homeq \Fun_x(\mathsf{C},T(\spaces)) \]
the $\infty$-category of functors which lift $x$ to $T(\spaces)$ along the projection map $T(\spaces) \to \spaces$.
\end{example}

\begin{example}
Taking $\mathsf{C} = *$ in Example~\ref{ex:affine} we obtain an affine $\infty$-topos $\cat{A}^*$ whose $\infty$-category of points is $\spaces$. Therefore
\[ U(\cat{A}^*) \homeq \Fun^{\omega}(T(\spaces),\spaces) \]
and, for an $\infty$-groupoid $x \in \spaces$,
\[ p(U_x\cat{A}^*) \homeq T_x\spaces = \spectra(\spaces_{/x}) \]
the $\infty$-category of spectra parameterized over $x$. 
\end{example}


\bibliographystyle{amsalpha}
\bibliography{mcching}

\end{document}